\documentclass[12pt]{article}
\usepackage[amsmath]{e-jc}
\usepackage{graphicx}

\theoremstyle{remark}
\newtheorem{case}{Case}

\dateline{Jun 7, 2022}{Apr 18, 2023}{Jun 16, 2023}

\MSC{05C35, 05C65, 05D99, 82B43}

\Copyright{The authors. Released under the CC BY license (International 4.0).}
\title{On the Running Time of\\ Hypergraph Bootstrap Percolation}

\author{Jonathan A. Noel\thanks{Research supported by NSERC Discovery Grant RGPIN-2021-02460 and NSERC Early Career Supplement DGECR-2021-00024 and a Start-Up Grant from the University of Victoria.}\\
\small Department of Mathematics and Statistics\\[-0.8ex]
\small University of Victoria\\[-0.8ex] 
\small Victoria, B.C., Canada\\
\small\tt noelj@uvic.ca
\and
Arjun Ranganathan\thanks{Research supported by Kishore Vaigyanik Protsahan Yojana Fellowship (KVPY), Department of Science and Technology, Govt. of India. and by IISER Pune-IDeaS Scholarship}\\
\small Department of Mathematics\\[-0.8ex]
\small Indian Institute of Science Education and Research (IISER), Pune\\[-0.8ex]
\small Pune, India\\
\small\tt arjun.ranganathan@students.iiserpune.ac.in}

\begin{document}

\maketitle

% ABSTRACT
% E-JC papers must include an abstract. The abstract should consist of a
% succinct statement of background followed by a listing of the
% principal new results that are to be found in the paper. The abstract
% should be informative, clear, and as complete as possible. Phrases
% like "we investigate..." or "we study..." should be kept to a minimum
% in favor of "we prove that..."  or "we show that...".  Do not
% include equation numbers, unexpanded citations (such as "[23]"), or
% any other references to things in the paper that are not defined in
% the abstract. The abstract may be distributed without the rest of the
% paper so it must be entirely self-contained.  Try to include all words
% and phrases that someone might search for when looking for your paper.

\begin{abstract}
  Given $r\geq2$ and an $r$-uniform hypergraph $F$, the \emph{$F$-bootstrap process} starts with an $r$-uniform hypergraph $H$ and, in each time step, every hyperedge which ``completes'' a copy of $F$ is added to $H$. The maximum running time of this process has been recently studied in the case that $r=2$ and $F$ is a complete graph by Bollob\'as, Przykucki, Riordan and Sahasrabudhe~[Electron. J. Combin. 24(2) (2017), Paper No. 2.16], Matzke~[arXiv:1510.06156v2] and Balogh, Kronenberg, Pokrovskiy and Szab\'o~[arXiv:1907.04559v1]. We consider the case that $r\geq3$ and $F$ is the complete $r$-uniform hypergraph on $k$ vertices. Our main results are that the maximum running time is $\Theta\left(n^r\right)$ if $k\geq r+2$ and $\Omega\left(n^{r-1}\right)$ if $k=r+1$. For the case $k=r+1$, we conjecture that our lower bound is optimal up to a constant factor when $r=3$, but suspect that it can be improved by more than a constant factor for large $r$.
\end{abstract}

\section{Introduction}

Bootstrap percolation is a model of propagation phenomena in discrete structures which originated in the statistical physics literature to study the dynamics of ferromagnetism~\cite{ChalupaLeathReich79,AdlerStaufferAharony89,AizenmanLebowitz88,vanEnter87}. It has since inspired the introduction of many related models that have been the subject of intense research in sociology~\cite{Granovetter78,Watts02}, computer science~\cite{Flocchini04,Dreyer09} and, especially, mathematics~\cite{Bollobas+17,BaloghBollobasMorris12,HambardzumyanHatamiQian20,BenevidesPrzykucki15,BenevidesPrzykucki13,Przykucki12,Hartarsky18,BalisterBollobasSmith16,MorrisonNoel21,MorrisonNoelScott17,BaloghBollobasMorrisRiordan12,BollobasSmithUzzell15,BaloghBollobas06,MorrisonNoel18,GravnerHolroyd08,Holroyd03,GravnerHolroydSivakoff21,Balogh+12,ShapiroStephens91,Holroyd06,BaloghBollobasMorris09-maj,CerfCirillo99,HartarskyMorris19,Kolesnik22,CokerGunderson14,GundersonKochPrzykucki17}. 

Our focus in this paper is on \emph{hypergraph bootstrap percolation} which generalizes the notion of \emph{graph bootstrap percolation} introduced by Bollob\'as~\cite{Bollobas68} in 1968, originally under the name ``weak saturation.'' For a set $S$ and $r\geq0$, let $S^{\left(r\right)}$ denote the collection of all subsets of $S$ of cardinality $r$.  Given an $r$-uniform hypergraph $F$, the \emph{$F$-bootstrap process} starts with an $r$-uniform hypergraph $H_0$. Then, at each time step $t\geq1$, the hypergraph $H_t$ consists of all hyperedges of $H_{t-1}$ as well as each $e\in V\left(H_0\right)^{\left(r\right)}$ for which there exists a copy $F'$ of $F$ containing $e$ such that $E\left(F'\right)\subseteq E\left(H_{t-1}\right)\cup \left\{e\right\}$ (that is, $e$ ``completes'' a copy of $F$ when added to $H_{t-1}$). The hypergraph $H_0$ is referred to as the \emph{initial infection} and, for $t\geq0$, the hyperedges of $H_t$ are said to be \emph{infected at time $t$}. 

The most well studied problem related to the $F$-bootstrap process is to determine the minimum number of hyperedges in an initial infection $H_0$ on $n$ vertices such that every element of $V\left(H_0\right)^{\left(r\right)}$ is eventually infected. When $F$ is the complete $r$-uniform hypergraph on $k$ vertices, which we denote by $K_k^r$, this is equivalent to the famous Skew Two Families Theorem proved by Frankl~\cite{Frankl82} using a beautiful exterior algebraic approach of Lov\'asz~\cite{Lovasz77}; see also~\cite{Kalai85,AlonKalai85}. Alon~\cite{Alon85} proved a multipartite version of the Skew Two Families Theorem, which was extended further by Moshkovitz and Shapira~\cite{MoshkovitzShapira15}; the results of both of these papers can be thought of in terms of hypergraph bootstrap percolation in a multipartite ``host'' hypergraph. Pikhurko~\cite{Pikhurko01,Pikhurko01a} solved the problem for several classes of non-complete hypergraphs $F$ using the exterior algebraic method of~\cite{Lovasz77,Frankl82} and a related matroid-theoretic approach of Kalai~\cite{Kalai85}. For results on graphs, i.e. the case $r=2$, see~\cite{Bollobas68,KronenbergMartinsMorrison21,Bidgoli+21,MorrisonNoelScott17,Kalai84,Kalai85}.

Another major focus in bootstrap percolation has been on analyzing the maximum number of time steps that a given bootstrap percolation model can take before it stabilizes~\cite{Bollobas+17,Matzke15,BenevidesPrzykucki15,BenevidesPrzykucki13,Przykucki12,Hartarsky18,BalisterBollobasSmith16}. In the context of the $F$-bootstrap process  with initial infection $H_0$, the \emph{running time} is
\[M_F\left(H_0\right):=\min\left\{t: H_t=H_{t+1}\right\}\]
and the \emph{maximum running time} among all hypergraphs $H_0$ with $n$ vertices is denoted by $M_F\left(n\right)$. If $F=K_k^r$ for some $k$, then we write $M_F\left(H_0\right)$ and $M_F\left(n\right)$ as $M_k^r\left(H_0\right)$ and $M_k^r\left(n\right)$, respectively. We omit the superscript in the case $r=2$. 

Our goal is to bound $M_k^r\left(n\right)$ for fixed $r$ and $k$ asymptotically as a function of $n$. This problem was first studied in the case $r=2$ independently by Bollob\'as, Przykucki, Riordan and Sahasrabudhe~\cite{Bollobas+17} and Matzke~\cite{Matzke15}. Clearly, $M_k^r\left(n\right)\leq \binom{n}{r}$ for all $n,r$ and $k$ as at least one hyperedge must become infected in each step. It is an easy exercise to show that the maximum running time for the $K_3$-bootstrap process, i.e. $M_3\left(n\right)$, is precisely $\lceil \log_2\left(n-1\right)\rceil$. Interestingly, for $k=4$, the maximum running time jumps from logarithmic to linear.

\begin{theorem}[Bollob\'as et al.~\cite{Bollobas+17}, Matzke~\cite{Matzke15}]
\label{th:K4}
$M_4\left(n\right)=n-3$ for all $n\geq3$.
\end{theorem}

Bollob\'as et al.~\cite{Bollobas+17} conjectured that $M_k\left(n\right)=o\left(n^2\right)$ for all $k\geq5$. This was disproved for all $k\geq6$ by Balogh, Kronenberg, Pokrovskiy and Szab\'o~\cite{Balogh+19}. Since the inequality $M_{k+1}\left(n+1\right)\geq M_k\left(n\right)$ can be shown by simply taking a construction on $n$ vertices, adding a vertex and infecting every edge incident to it (see~\cite[Proposition~10]{Balogh+19} or Lemma~\ref{lem:easyStepUp} in this paper), the critical case is $k=6$. 

\begin{theorem}[Balogh et al.~\cite{Balogh+19}]
\label{th:K6}
$M_6\left(n\right)\geq \frac{n^2}{2500}$ for all $n$ sufficiently large.
\end{theorem}

The growth rate of the maximum running time for the $K_5$-bootstrap process is still unknown. In fact, it is an interesting open problem to determine whether or not it is quadratic~\cite{Bollobas+17,Balogh+19}. To date, the best known lower bound is given by Balogh et al.~\cite{Balogh+19} by exploiting connections to additive combinatorics; their result improved on a bound of $n^{13/8-o\left(1\right)}$ due to~\cite{Bollobas+17}. Let $r_3\left(n\right)$ be the largest cardinality of a subset of $[n]$ without a $3$-term arithmetic progression. While the asymptotics of $r_3\left(n\right)$ are not known precisely, Roth's Theorem~\cite{Roth53} implies that $r_3\left(n\right)=o\left(n\right)$ (see~\cite{Schoen21,KellyMeka23+} for recent quantitative bounds) and the Behrend Construction~\cite{Behrend} yields $r_3\left(n\right)\geq n^{1-O\left(1/\sqrt{\log\left(n\right)}\right)}$.

\begin{theorem}[Balogh et al.~\cite{Balogh+19}]
\label{th:K5}
$M_5\left(n\right)\geq \frac{nr_3\left(n\right)}{1200}$. 
\end{theorem}

Here, we initiate the study of $M_k^r\left(n\right)$ for $r\geq3$. Our main contributions are constructions of initially infected hypergraphs yielding lower bounds. Our first result concerns the case that $k=r+1$. It would be interesting to know whether this bound is tight up to a constant factor, especially in the case $r=3$ and $k=4$; see Conjecture~\ref{conj:K43} and Question~\ref{ques:gap1}.

\begin{theorem}
\label{th:n^2}
Let $r\geq3$. If $k=r+1$, then $M_k^r\left(n\right)=\Omega\left(n^{r-1}\right)$. 
\end{theorem}

In contrast, for $k\geq r+2$, we show that the trivial upper bound $\binom{n}{r}$ is tight up to a constant factor. 

\begin{theorem}
\label{th:n^3}
Let $r\geq3$. If $k\geq r+2$, then $M_k^r\left(n\right)=\Theta\left(n^{r}\right)$. 
\end{theorem}

The rest of the paper is organized as follows. In the next section, we build up some basic notation and terminology and establish a few preliminary lemmas. In particular, these lemmas will be used to reduce Theorems~\ref{th:n^2} and~\ref{th:n^3} to finding constructions with some additional properties for $r=3$ and $k\in\left\{4,5\right\}$. In Section~\ref{sec:beachball}, we introduce a key construction, which we call the ``beachball hypergraph''. This hypergraph has running time that is only linear with respect to its number of vertices, but plays a key role in the proofs of both of our main theorems. In the same section, we show how linearly many beachball hypergraphs can be ``chained together'' to prove Theorem~\ref{th:n^2}. In Section~\ref{sec:cubic}, we use the beachball construction in a different way to prove Theorem~\ref{th:n^3}. We conclude the paper in Section~\ref{sec:concl} with two open problems.

\section{Preliminaries}
\label{sec:prelim}

Given $r$-uniform hypergraphs $F$ and $H$, say that a copy $F'$ of $F$ is \emph{susceptible} to $H$ if there exists a hyperedge $e\notin E\left(H\right)$  and $e\in E\left(F'\right)$ such that $E\left(F'\right)\subseteq E\left(H\right)\cup \left\{e\right\}$. Say that $H$ is \emph{$F$-stable} if there are no copies of $F$ that are susceptible to $H$. We make the following simple observation.

\begin{observation}
Let $H_0$ be an initial infection for the $F$-bootstrap process. Then, for $t\geq1$, $E\left(H_t\right)$ is the union of $E\left(H_{t-1}\right)$ and the hyperedge sets of all copies of $F$ which are susceptible to $H_{t-1}$. 
\end{observation}

The following very straightforward definition and the lemma that follows it allow us to use a lower bound on $M_k^r\left(n\right)$ to get a lower bound on $M_{k+1}^{r}\left(n\right)$ which is only slightly worse; c.f.~\cite[Proposition~10]{Balogh+19}. 

\begin{definition}
Given an $r$-uniform hypergraph $H$ and $w\notin V\left(H\right)$, let $H\vee w$ be the hypergraph with vertex set $V\left(H\right)\cup \left\{w\right\}$ and hyperedge set
\[E\left(H\right)\cup \left\{e\cup\left\{w\right\}: e\in V\left(H\right)^{\left(r-1\right)}\right\}.\]
\end{definition}

\begin{lemma}
\label{lem:easyStepUp}
Let $H_0$ be an initial infection for the $K_k^r$-bootstrap process and, for $w\notin V\left(H_0\right)$, let $H_0'=H_0\vee w$ be an initial infection for the $K_{k+1}^r$-bootstrap process. Then  $H_t'=H_t\vee w$ for all $t\geq0$. 
\end{lemma}

\begin{proof}
Suppose not and let $t$ be the minimum time that the equality is violated. Since $H_0'=H_0\vee w$ by definition, we must have $t\geq 1$. Observe that no hyperedge containing $w$ becomes infected in any step $t\geq1$ since all such hyperedges are already in $H_0$.

First, suppose that there is $e\in E\left(H_t \vee w\right)$ such that $e\notin E\left(H_t'\right)$. As noted above, we may assume that $e \in E\left(H_t\right)$ since $w\notin e$. Let $F$ be the corresponding copy of $K_k^r$ that is susceptible to $H_{t-1}$. Then, by minimality of $t$, all hyperedges of $F\vee w$ except for $e$ are present in $H_{t-1}'$. Thus, $e\in E\left(H_t'\right)$, which is a contradiction. The proof of the other direction is similar.
\end{proof}

The analysis of the running time of $H_0$ tends to become unwieldy if there are time steps $t$ in which there is more than one copy of $F$ that is susceptible to $H_t$. For this reason, the constructions in this paper will be designed with a specific goal of avoiding this situation; a similar approach is taken for graphs in~\cite{Balogh+19,Bollobas+17}. This motivates the following definitions.

\begin{definition}
Let $H_0$ be an initial infection for the $F$-bootstrap process. Say that $H_0$ is \emph{$F$-tame} if there is at most one copy of $F$ which is susceptible to $H_t$ for all $t\geq0$. 
\end{definition}

\begin{definition}
Let $H_0$ be $F$-tame. The \emph{trajectory} of $H_0$ is the sequence
\[\left(F_0,e_1,F_1,\dots,e_{T-1},F_{T-1},e_T\right)\]
where $T=M_F\left(H_0\right)$ and, for all $0\leq t\leq T-1$, $F_t$ is the unique copy of $F$ that is susceptible to $H_t$ and $e_{t+1}$ is the unique hyperedge of $E\left(F_t\right)\setminus E\left(H_t\right)$. 
\end{definition}

The following proposition is an easy consequence of Lemma~\ref{lem:easyStepUp}; we omit the proof.

\begin{proposition}
\label{prop:easyStepUp}
If $H_0$ is $K_k^r$-tame with trajectory $(F_0,e_1,F_1,\dots,F_{T-1},e_T)$, then $H_0\vee w$ is $K_{k+1}^{r}$-tame with trajectory $(F_0\vee w, e_1,F_1\vee w,\dots,F_{T-1}\vee w, e_T)$.
\end{proposition}

In our constructions, it will be especially useful to restrict our attention to $F$-tame hypergraphs with limited interaction between the elements of their trajectories.  

\begin{definition}
\label{defn:civilized}
Let $\left(F_0,e_1,\dots,F_{T-1},e_T\right)$ be the trajectory of an $F$-tame hypergraph $H_0$ and let $e_0\in E\left(F_0\right)$. We say that $H_0$ is \emph{$F$-civilized} with respect to $e_0$ if the following two conditions hold:
\begin{enumerate}
    \item[(a)] \label{civilized1} $E\left(F_j\right)\cap \left\{e_0,e_1,\dots,e_T\right\}=\left\{e_j,e_{j+1}\right\}$ for any $0\leq j\leq T-1$ and
    \item[(b)] \label{civilized2} $H_0\setminus\left\{e_0\right\}$ is $F$-stable. 
\end{enumerate}
\end{definition}

In other words, an $F$-tame hypergraph $H_0$ is $F$-civilized with respect to $e_0$ if every copy of $F$ in its trajectory has exactly two hyperedges missing from $H_0\setminus\left\{e_0\right\}$ and $H_0\setminus\left\{e_0\right\}$ is $F$-stable. A useful property of an $F$-civilized hypergraph is that its trajectory can be reversed by swapping $e_0$ and $e_T$.

\begin{lemma}
\label{lem:reverse}
Let $H_0$ be a hypergraph which is $F$-civilized with respect to $e_0$ with trajectory $\left(F_0,e_1,\dots,F_{T-1},e_T\right)$ and let $H_0'=H_0\setminus\left\{e_0\right\}\cup \left\{e_T\right\}$. Then $H_0'$ is $F$-civilized with respect to $e_T$ with trajectory $\left(F_{T-1},e_{T-1},\dots,F_0,e_0\right)$.
\end{lemma}

\begin{proof}

To show that $H_0'$ is $F$-civilized with respect to $e_T$, it suffices to show that it is $F$-tame with the correct trajectory since the additional conditions of Definition~\ref{defn:civilized} will follow from the fact that $H_0$ is $F$-civilized with respect to $e_0$. 

For $0\leq t< T$, let $F'$ be a copy of $F$ which is susceptible to $H_{t}'$ and let $e'$ be the unique hyperedge of $F'$ that is not in $H_{t}'$. We claim that $F'=F_{T-t-1}$ and $e'=e_{T-t-1}$. For the sake of contradiction, suppose that at least one of these equalities does not hold and let $t$ be the minimum time for which such an $F'$ and $e'$ exist. By minimality of $t$, we have that
\[H_t' = H_0'\cup\left\{e_{T-t}, \dots,e_{T-1}\right\}.\]
In particular, $H_t'$ is a subhypergraph of $H_T$. Now, observe that $e'$ must be contained in $H_T$; if not, then $F'$ would be susceptible to $H_T$ which contradicts the assumption that $\left(F_0,e_1,\dots,F_{T-1},e_T\right)$ is the trajectory of $H_0$. Next, we claim that $F'$ is not a subhypergraph of $H_0$. If it were, then, since $e'$ is not contained in $H_t'$, we must have that $e'=e_0$ and all other hyperedges of $F'$ are contained in $H_0$. However, this would contradict condition \ref{civilized2} of Definition~\ref{defn:civilized}. Therefore, $F'$ is not a subhypergraph of $H_0$.
%\textcolor{red}{\left(I think we should say that $t$ is the first time the statement of the theorem fails\right)} 
%\textcolor{red}{\left(Yes, you are right.\right)} 

So, $H_0$ does not contain every hyperedge of $F'$ but $H_T$ does. Thus, we can let $0\leq j<T$ be the maximum index such that there is a hyperedge of $F'$ that is not contained in $H_j$. So, $F'$ is susceptible to $H_j$, which, by the assumption of the lemma, implies that $F'=F_j$. Since $H_0$ is $F$-civilized with respect to $e_0$, the only hyperedges of $F_j$ missing from $H_0\setminus\left\{e_0\right\}$ are $e_j$ and $e_{j+1}$. Since $F_j$ is susceptible to $H_{t}'=H_0'\cup\left\{e_{T-t},\dots,e_{T-1}\right\}$ and also to $H_j$, but not to $H_{t-1}'$, the only possibility is that $e_{T-t}=e_{j+1}$ and $e'=e_j$. 
%\textcolor{red}{\left(assuming minimality of $t$ again right?\right)}. 
%\textcolor{red}{\left(Yep, good catch!\right)}. 
So, $j+1=T-t$ which implies that $F'=F_j=F_{T-t-1}$ and $e_j=e_{T-t-1}$, as we wanted. 

Thus, $F_{T-t-1}$ is the only copy of $F$ that can be susceptible to $H_t'$. The last thing to show is that $F_{T-t-1}$ is, indeed, susceptible to $H_t'$. To see this, note that $H_0$ is $F$-civilized with respect to $e_0$, and so $e_{T-t-1}$ is the unique hyperedge of $F_{T-t-1}$ that is not in $H_t'$. Thus, $H_0'$ is $F$-tame with the desired trajectory and the proof is complete. 
\end{proof}

Next, we define a natural ``step up'' construction for converting constructions for the $K_k^r$-bootstrap process into constructions for the $K_{k+1}^{r+1}$-bootstrap process. Later, in Lemma~\ref{lem:stepUp}, we will see how Lemma~\ref{lem:reverse} can be applied to ``chain together'' these step up constructions to transform a $K_k^r$-civilized hypergraph into a $K_{k+1}^{r+1}$-civilized hypergraph with a linear number of extra vertices in such a way that the running time is boosted by a linear factor. 

\begin{definition}
Let $r\geq2$, let $H$ be an $r$-uniform hypergraph and let $w\notin V\left(H\right)$. Define $H^{+w}$ to be the $\left(r+1\right)$-uniform hypergraph with vertex set $V\left(H\right)\cup\left\{w\right\}$ and hyperedges
\[V\left(H\right)^{\left(r+1\right)}\cup\left\{e\cup\left\{w\right\}: e\in E\left(H\right)\right\}.\]
\end{definition}

Next, we show that 
%$H_0^{\vee w}$ and \textcolor{red}{This one hasn't been defined yet.}
$H_0^{+w}$ 
%behave 
behaves in essentially the same way with respect to the
%$K_{k+1}^{r}$ and 
$K_{k+1}^{r+1}$-bootstrap 
%processes respectively 
process as $H_0$ does with respect to the $K_k^r$-bootstrap process.

\begin{lemma}
\label{lem:simpleStepUp}
For $k\geq r\geq2$, let $H_0$ be an initial infection for the $K_k^r$-bootstrap process, let $w\notin V\left(H_0\right)$, and let $H_0' = H_0^{+w}$ be the initial infection for the $K_{k+1}^{r+1}$-bootstrap process. Then $H_t'=\left(H_t\right)^{+w}$ for all $t\geq0$. 

\end{lemma}

\begin{proof}
Suppose that the lemma is false and let $t$ be the smallest index such that $H_t'\neq \left(H_t\right)^{+w}$. By definition, we have $H_0'=H_0^{+w}$ and so $t\geq 1$. 

First, suppose that there is a hyperedge which is contained in $\left(H_t\right)^{+w}$ but not in $H_t'$. Since both of these hypergraphs contain all of $V\left(H_0\right)^{\left(r+1\right)}$, this hyperedge must have the form $e\cup\left\{w\right\}$ where $e\in E\left(H_t\right)$. By minimality of $t$, we must have $e\notin E\left(H_{t-1}\right)$ and so there must be a copy $F$ of $K_k^r$ containing $e$ which is susceptible to $H_{t-1}$. As $F$ is a copy of $K_k^{r}$ containing $e$, this means that $F^{+w}$ must be a copy of $K_{k+1}^{r+1}$ containing $e\cup\left\{w\right\}$. However, by minimality of $t$, we have that $H_{t-1}'=\left(H_{t-1}\right)^{+w}$, and hence $e\cup\left\{w\right\} \notin H_{t-1}'$. Since $H_{t-1}$ contains every hyperedge of $F$ except $e$, we conclude that $H_{t-1}'$ must contain every hyperedge of $F^{+w}$ except for $e\cup \left\{w\right\}$. Therefore, $F^{+w}$ is susceptible to $H_{t-1}'$ and so $e\cup\left\{w\right\}$ is contained in $H_t'$, which is a contradiction.

Now, suppose that there is a hyperedge $e'$ that is contained in $H_t'$ but not in $\left(H_t\right)^{+w}$. As in the previous case, this hyperedge must contain $w$. By minimality of $t$, we must have $e'\notin E\left(H_{t-1}'\right)$ and so there must be a copy $F'$ of $K_{k+1}^{r+1}$ containing $e'$ which is susceptible to $H_{t-1}'$. However, by minimality of $t$, we have $H_{t-1}'=\left(H_{t-1}\right)^{+w}$ and so $H_{t-1}$ contains all hyperedges of the form $e''\setminus\left\{w\right\}$ for $e''\in E\left(F'\right)\setminus\left\{e'\right\}$. Thus, the copy of $K_k^r$ on vertex set $V\left(F\right)\setminus\left\{w\right\}$ is susceptible to $H_{t-1}$ which implies that $e'$ is in $\left(H_t\right)^{+w}$. This contradiction completes the proof. 
\end{proof}

As alluded to earlier, we now apply Lemmas~\ref{lem:reverse} and~\ref{lem:simpleStepUp} to show that step up constructions for $K_k^r$-civilized hypergraphs can be chained together to yield $K_{k+1}^{r+1}$-civilized hypergraph with very long running time. 

\begin{lemma}
\label{lem:stepUp}
For $r\geq3$ and $k\geq r+1$, let $H_0$ be a $K_k^r$-civilized hypergraph. Then, for any $m\geq1$, there exists a $K_{k+1}^{r+1}$-civilized hypergraph $H_0'$ with $|V\left(H_0\right)|+m + \left(m-1\right)\left(k-r-1\right)$ vertices such that
\[M_{k+1}^{r+1}\left(H_0'\right)=m\cdot M_k^r\left(H_0\right) + m-1.\]
\end{lemma}

\begin{proof}
Let $e_0$ be a hyperedge of $H_0$ such that $H_0$ is $K_k^r$-civilized with respect to $e_0$ and let $\left(F_0,e_1,\dots,F_{T-1},e_T\right)$ be the trajectory of $H_0$. Let us describe the construction of $H_0'$. First, let $w_1,\dots,w_m$ be distinct vertices which are not in $V\left(H_0\right)$. For any $e\in V\left(H_0\right)^{\left(r\right)}$ and $1\leq j\leq m$, let $e^{+j}=e\cup\left\{w_j\right\}$. Let $H_0^1$ be $H_0^{+w_1}$ and, for $2\leq j\leq m$, let $H_0^j$ be $H_0^{+w_j}\setminus\left\{e_0^{+j}\right\}$. For each $1\leq j\leq m-1$, let $X_j:=\left\{x_j^1,\dots,x_j^{k-r-1}\right\}$ be a set of vertices disjoint from the set of vertices introduced so far. For $1\leq j\leq m-1$, let $f\left(j\right)=T$ if $j$ is odd and $f\left(j\right)=0$ if $j$ is even. For $1\leq j\leq m-1$, let $C_j$ be the $\left(r+1\right)$-uniform hypergraph with vertex set $\left\{w_j,w_{j+1}\right\}\cup X_j\cup e_{f\left(j\right)}$ containing all hyperedges in this set, except for $e_{f\left(j\right)}^{+j}$ and $e_{f\left(j\right)}^{+\left(j+1\right)}$. Finally, 
\[H_0':=\left(\bigcup_{j=1}^mH_0^{j}\right)\cup\left(\bigcup_{j=1}^{m-1}C_j\right).\] 
Intuitively, the way to think of this is as follows. For each $1\leq j\leq m$, the hypergraph $H_0^j$ will emulate the $K_k^{r}$-bootstrap process with initial infection $H_0$, either in the forward or backwards (as in Lemma~\ref{lem:reverse}) direction, depending on the parity of $j$. The hypergraph $C_j$ for $1\leq j\leq m-1$ is used to link the process on $H_0^j$ to the process $H_0^{j+1}$ in such a way that the termination of the former triggers the start of the latter. 

Let us formalize this. First, we show that, if $x\in X_j$ for some $1\leq j\leq m-1$, then, for all $t\geq0$, every hyperedge $e$ of $H_t'$ containing $x$ is contained in $V\left(C_j\right)$. Suppose that this is not the case, let $t$ be the minimum time for which it fails and let $z$ be a vertex of $e$ that is not in $V\left(C_j\right)$. Clearly, by construction of $H_0'$, we must have $t\geq1$. By minimality of $t$, there must be a copy $F$ of $K_{k+1}^{r+1}$ which is susceptible to $H_{t-1}'$ and contains $e$. Since $k+1 > r+1$, we can let $u$ be a vertex of $V\left(F\right)\setminus e$. Since $r+1\geq4$, we can choose $y$ to be a vertex of $e\setminus\left\{x,z\right\}$. Now, consider the set $e'=e\setminus\left\{y\right\}\cup\left\{u\right\}$. This is a hyperedge of $F$ which contains both of $x$ and $z$. However, all hyperedges of $F$ other than $e$ are in $H_{t-1}'$, and so $e'$ is in $H_{t-1}'$ which contradicts the minimality of $t$. 

Next, observe that, for any $t\geq0$, the hypergraph $H_t'$ does not contain any hyperedge which includes $w_j$ and $w_{j'}$ for $|j-j'|>1$. Note that $H_0'$ has no such hyperedge. Consider the first time $t\geq1$ that such a hyperedge appears, and observe that the relevant copy of $K_{k+1}^{r+1}$ that was allegedly susceptible to $H_{t-1}'$ is missing all other hyperedges containing $w_j$ and $w_{j'}$ by minimality of $t$, which is a contradiction. 

Now, let us show that, for $1\leq j\leq m-1$, for every $t\geq0$, every hyperedge $e$ of $H_t'$ containing $w_j$ and $w_{j+1}$ is contained in $V\left(C_j\right)$. If not, let $t$ be the minimum time that it is violated and let $z$ be the offending vertex of $e$. As usual, $t\geq1$ by construction of $H_0'$. Take a copy $F$ of $K_{k+1}^{r+1}$ which is susceptible to $H_{t-1}'$ and contains $e$. Since $k+1>r+1$, we can let $u$ be a vertex of $V\left(F\right)\setminus e$. Since $r+1\geq4$, we can pick $y\in e\setminus\left\{w_j,w_{j+1},z\right\}$. Now, consider $e'=e\setminus \left\{y\right\}\cup\left\{u\right\}$. This is a hyperedge of $F$ other than $e$ which contains $w_j,w_{j+1}$ and $z$, which contradicts the minimality of $t$. 

Putting all of this together, we see that the only hyperedges that can become infected during the $K_{k+1}^{r+1}$-bootstrap process starting with $H_0'$ are those which are contained within $V\left(H_0^j\right)$ for some $1\leq j\leq m$ or within $V\left(C_j\right)$ for some $1\leq j\leq m-1$. The hyperedges of $H_0'$ within $V\left(C_j\right)$ form nothing more than a copy of $K_{k+1}^{r+1}$ with two hyperedges removed. The same is true for the hyperedges within $V\left(F_i\right)\cup\left\{w_j\right\}$ for any copy $F_i$ of $K_k^r$ in the trajectory of $H_0$ and any $1\leq j\leq m$ with the exception of the case $i=0$ and $j=1$, in which case the hyperedges form a copy of $K_{k+1}^{r+1}$ with one hyperedge removed. Therefore, the trajectory of $H_0'$ follows that of $H_0^{+w_1}$ for the first $M_k^r\left(H_0\right)$ steps. At that point, the hyperedge $e_T^{+1}$ has become infected, which causes $C_1$ to become susceptible. Thus, $e_T^{+2}$ becomes infected. This triggers the process in $H_0^{2}$ to start running in the reverse direction as in Lemma~\ref{lem:reverse}, and so on. Putting all of this together, the fact that $H_0$ is $K_k^r$-civilized with respect to $e_0$ now implies that $H_0'$ is $K_{k+1}^{r+1}$-civilized with respect to $e_0^{+1}$. The running time includes $M_k^r\left(H_0\right)$ steps for every $1\leq j\leq m$, plus one extra step for each $1\leq j\leq m-1$. This completes the proof. 
\end{proof}

Using the results of this section, we show that, in order to prove Theorems~\ref{th:n^2} and~\ref{th:n^3}, it suffices to find $K_k^r$-civilized constructions with long running time for $r=3$ and $k\in \left\{4,5\right\}$. 

\begin{corollary}
\label{cor:reduction}
If there is a $K_4^3$-civilized hypergraph $H_0$ on $n$ vertices with $M_4^3\left(H_0\right)=\Omega\left(n^2\right)$, then Theorem~\ref{th:n^2} holds. Likewise, if there is a $K_5^3$-civilized hypergraph $H_0'$ on $n$ vertices with $M_5^3\left(H_0'\right)=\Omega\left(n^3\right)$, then Theorem~\ref{th:n^3} holds.
\end{corollary}

\begin{proof}
For any fixed $r\geq3$ and $k$ such that $k=r+1$, applying Lemma~\ref{lem:stepUp} with $m=\Theta\left(n\right)$ exactly $r-3$ times to the hypergraph $H_0$ in the statement of the corollary yields an $r$-uniform hypergraph with $O\left(n\right)$ vertices whose running time with respect to the $K_k^r$-bootstrap process is $\Omega\left(n^{r-1}\right)$. For $r\geq3$ and $k\geq r+2$, apply Proposition~\ref{prop:easyStepUp} to $H_0'$ exactly $k-r-2$ times, and then Lemma~\ref{lem:stepUp} with $m=\Theta\left(n\right)$ exactly $r-3$ times to get an $r$-uniform hypergraph with $O\left(n\right)$ vertices whose running time with respect to the $K_k^r$-bootstrap process is $\Omega\left(n^r\right)$.  Trivially, $M_k^r\left(n\right)\leq \binom{n}{r}=O\left(n^r\right)$ and so $M_k^r\left(n\right)=\Theta\left(n^r\right)$. 
\end{proof}

\section{Getting the Beachball Rolling}
\label{sec:beachball}

\begin{figure}[htbp]
  \centering
    \includegraphics[width=0.60\textwidth]{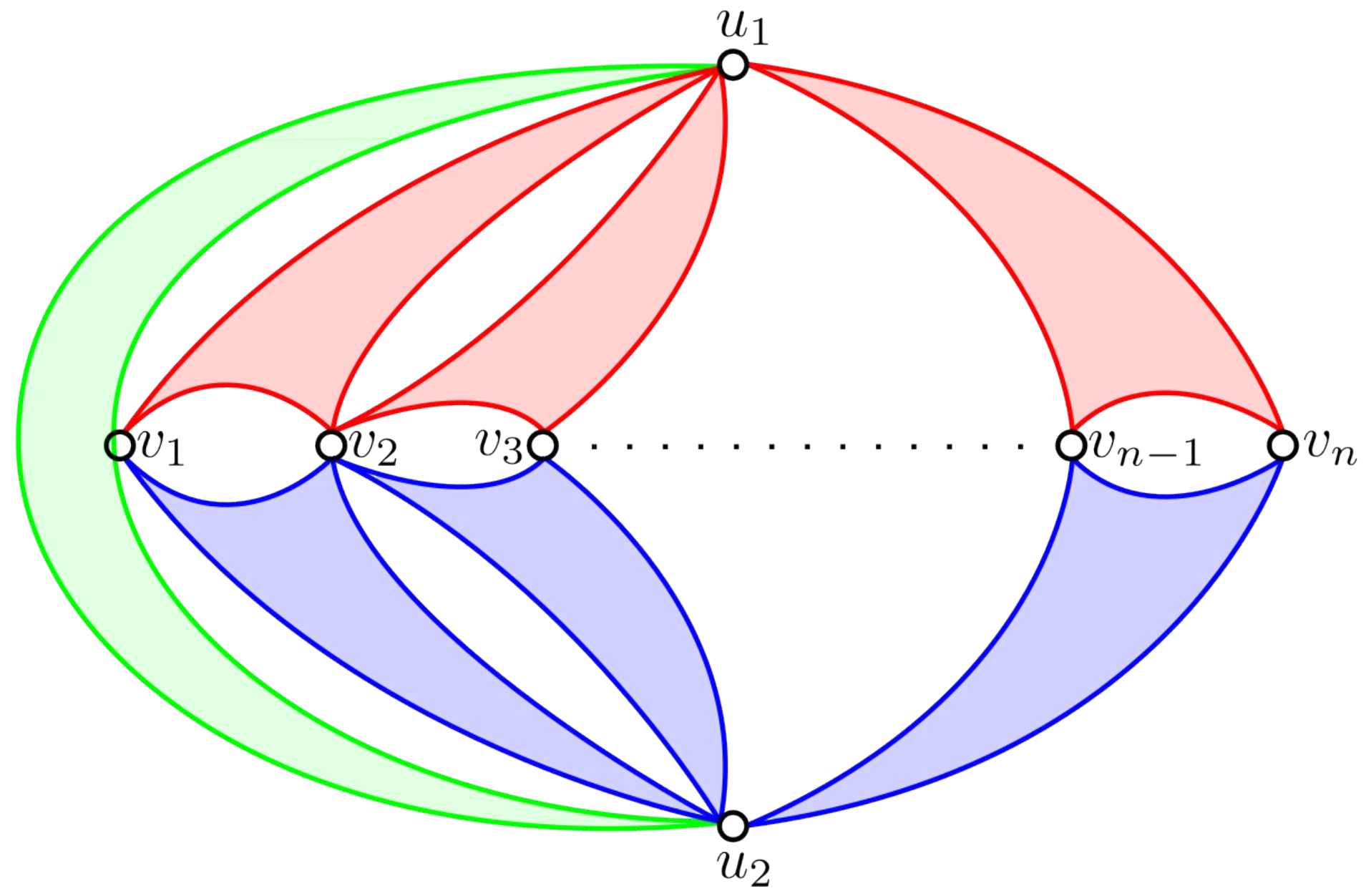}
    \caption{An illustration of the hypergraph $B\left(v_1,\dots,v_n,u_1,u_2\right)$.}
    \label{fig:beachball}
\end{figure}

The following definition provides a simple gadget that will be used in both of our main constructions. See Figure~\ref{fig:beachball} for an illustration. 

\begin{definition}
For $n\geq 1$ and distinct vertices $v_1,v_2,\dots,v_n,u_1$ and $u_2$, the \emph{beachball hypergraph}  $B\left(v_1,\dots,v_n,u_1,u_2\right)$ is the $3$-uniform hypergraph with vertex set $\left\{v_1,\dots,v_n\right\}\cup\left\{u_1,u_2\right\}$ and hyperedges $\left\{u_1,v_i,v_{i+1}\right\}$ and $\left\{u_2,v_i,v_{i+1}\right\}$ for $1\leq i\leq n-1$, as well as the hyperedge $\left\{u_1,u_2,v_1\right\}$. 
\end{definition}

We call $u_1$ and $u_2$ the \emph{top} and \emph{bottom} vertices of $B\left(v_1,\dots,v_n,u_1,u_2\right)$, respectively. The vertices $v_1,v_2,\dots,v_n$ are the \emph{middle} vertices. A key property of beachball hypergraphs is that they behave in a very straightforward way with respect to the $K_4^3$-process. 

\begin{lemma}
\label{lem:bb}
For $n\geq2$, the hypergraph $B\left(v_1,\dots,v_n,u_1,u_2\right)$ is $K_4^3$-civilized with respect to $e_0=\left\{u_1,u_2,v_1\right\}$ with trajectory $\left(F_0,e_1,\dots,F_{n-2},e_{n-1}\right)$ where
\[V\left(F_i\right) = \left\{u_1,u_2,v_{i+1},v_{i+2}\right\},\qquad e_{i+1}=\left\{u_1,u_2,v_{i+2}\right\}\]
for all $0\leq i\leq n-2$. 
\end{lemma}

\begin{proof}
For the purposes of this proof, denote $B\left(v_1,\dots,v_n,u_1,u_2\right)$ by $B^n_0$. Say that a hyperedge $e$ is \emph{flat} if it is contained in $\left\{v_1,\dots,v_{n}\right\}$ and \emph{wide} if it contains exactly one of $u_1$ or $u_2$ and is not an element of $E\left(B^n_0\right)$. 

We claim that $B_t^n$ contains no flat or wide hyperedges for any $t\geq0$. For $t=0$, this is true by definition of $B_0^n$. Now, consider the minimum $t$ for which the statement does not hold and let $e\in B_t^n$ be a hyperedge which is either flat or wide. Then there must exist a copy $F$ of $K_4^3$ with $e\in E\left(F\right)$ which is susceptible to $B_{t-1}^n$. Let $v$ be the unique vertex of $V\left(F\right)\setminus e$. If $e$ is flat, then we can write $e=\left\{v_i,v_j,v_k\right\}$ for $1\leq i< j<k\leq n$.  If $v\in \left\{v_1,\dots,v_{n}\right\}$, then $B_{t-1}^n$ must contain the hyperedge $e\setminus\left\{v_i\right\}\cup \left\{v\right\}$. However, this hyperedge is flat, which contradicts the minimality of $t$. Likewise, if $v\in\left\{u_1,u_2\right\}$, then $B^n_{t-1}$ contains the hyperedge $\left\{v_i,v_k,v\right\}$, which is wide, and so we get another contradiction. Similarly, if $e$ is wide, then, depending on whether $v\in \left\{v_1,\dots,v_{n}\right\}$ or $v\in \left\{u_1,u_2\right\}$, we find either a flat or a wide hyperedge in $B^n_{t-1}$, respectively. This proves the claim. 

Thus, the only hyperedges which can become infected after any number of steps are those of the form $\left\{u_1,u_2,v_i\right\}$ for some $i\in\left\{2,\dots,n\right\}$. This implies that the only copies of $K_4^3$ that can be susceptible at any time are those whose vertex sets have the form $\left\{u_1,u_2,v_i,v_{i+1}\right\}$ for some $i\in \left\{1,\dots,n-1\right\}$. 

Now, we claim that, for all $0\leq t\leq n-2$, a hyperedge of the form $\left\{u_1,u_2,v_i\right\}$ for $1\leq i\leq n$ is contained in $B^n_t$ if and only if $i\leq t+1$. This is clearly true for $t=0$ by construction. Consider the minimum $t\geq 1$ for which this statement is false and let $j>t+1$ such that $e=\left\{u_1,u_2,v_j\right\}\in B^n_t$. Then, by the result of the previous paragraph, the only possible copies of $K_4^3$ containing $e$ which could be susceptible at time $t-1$ are the ones with vertex set $\left\{u_1,u_2,v_{j-1},v_j\right\}$ or $\left\{u_1,u_2,v_j,v_{j+1}\right\}$. However, this implies that either $\left\{u_1,u_2,v_{j-1}\right\}$ or $\left\{u_1,u_2,v_{j+1}\right\}$ is in $E\left(B^n_{t-1}\right)$; in either case, this contradicts the minimality of $t$. 

Therefore, $B^n_0$ is $K_4^3$-tame with the trajectory described in the statement in the lemma. Let us argue that $B_0^n$ is  $K_4^3$-civilized with respect to $e_0=\left\{u_1,u_2,v_1\right\}$. Condition \ref{civilized1} of Definition~\ref{defn:civilized} follows from the construction of the beachball hypergraph, since all of the sets $\left\{u_1,u_2,v_{i+1},v_{i+2}\right\}$ for $0\leq i\leq n-2$ form a copy of $K_4^3$ with two hyperedges missing from $B_0^n\setminus\left\{e_0\right\}$. For the same reason, condition \ref{civilized2} also holds. This completes the proof.
\end{proof}

In our construction for the $K_4^3$-bootstrap process, we will also require the hypergraph $B'\left(v_1,\dots,v_n,u_1,u_2\right)$ obtained from $B\left(v_1,\dots,v_n,u_1,u_2\right)$ by adding the hyperedge $\left\{u_1,v_1,v_n\right\}$. We call this an \emph{augmented beachball hypergraph}. Next, we show that the addition of this hyperedge extends the running time by one step without disrupting any of the steps that came before it. In fact, this hypergraph is still $K_4^3$-civilized, albeit with respect to a different hyperedge $e_0'\neq e_0$.

\begin{figure}[htbp]
  \centering
    \includegraphics[width=0.6\textwidth]{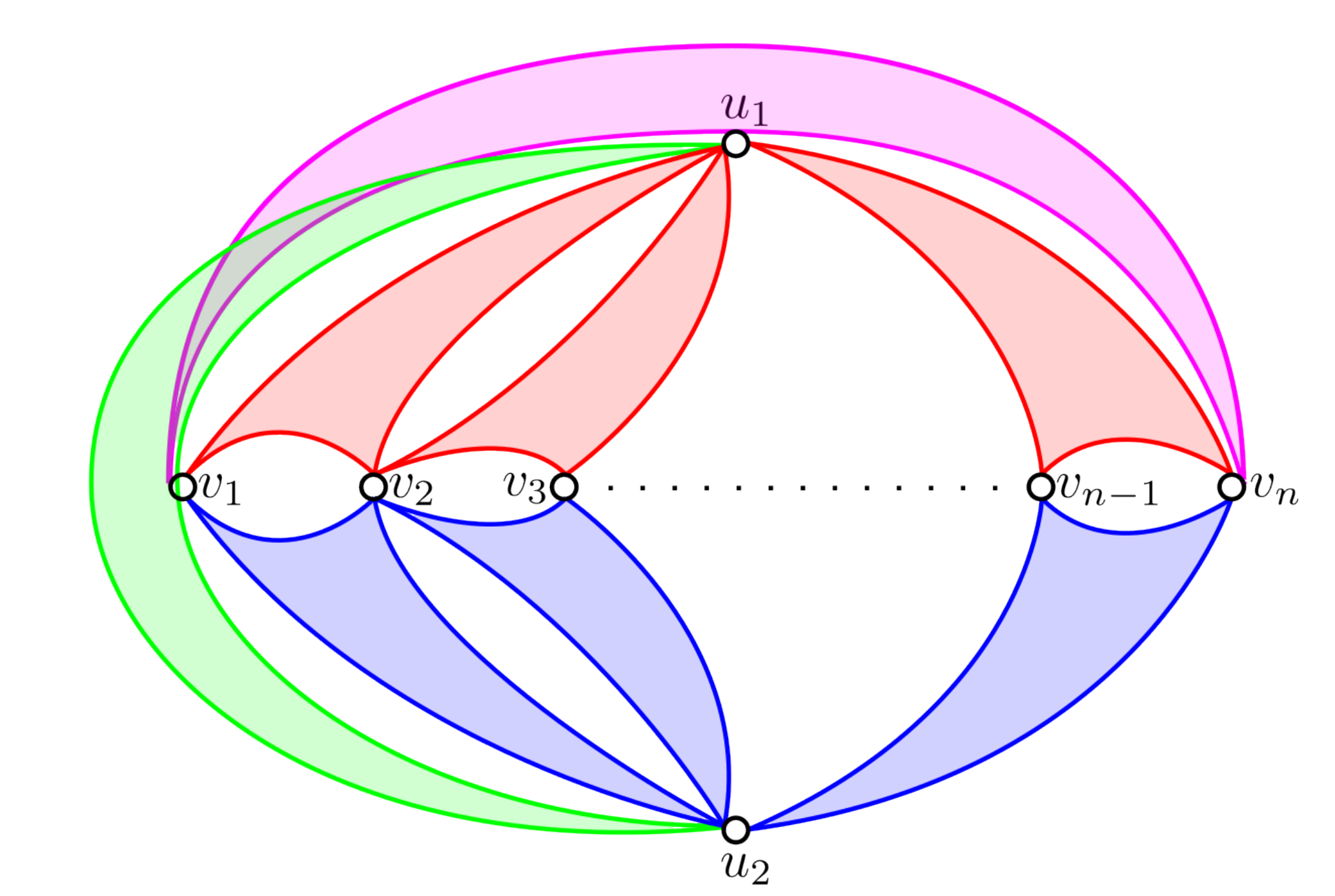}
    \caption{An illustration of the hypergraph $B'\left(v_1,\dots,v_n,u_1,u_2\right)$.}
    \label{fig:aug_beachball}
\end{figure}

\begin{lemma}
\label{lem:augbb}
For $n\geq4$, the hypergraph $B'\left(v_1,\dots,v_n,u_1,u_2\right)$ is $K_4^3$-civilized with respect to $e_0'=\left\{u_1,v_1,v_2\right\}$ with trajectory $\left(F_0,e_1,\dots,F_{n-1},e_{n}\right)$ where $\left(F_0,e_1,\dots,F_{n-2},e_{n-1}\right)$ is the trajectory of $B\left(v_1,\dots,v_n,u_1,u_2\right)$ and 
\[V\left(F_{n-1}\right) = \left\{u_1,u_2,v_1,v_n\right\},\qquad e_{n}=\left\{u_2,v_1,v_n\right\}.\]
\end{lemma}

\begin{proof}
Let us denote $B\left(v_1,\dots,v_n,u_1,u_2\right)$ by $B_0^n$ and $B'\left(v_1,\dots,v_n,u_1,u_2\right)$ by ${B'_0}^n$. By Lemma~\ref{lem:bb}, for every $0\leq t\leq n-2$, the only hyperedge of $B_t^n$ contained within $\left\{u_1,u_2,v_1,v_n\right\}$ is $\left\{u_1,u_2,v_1\right\}$. Therefore, if we additionally infect the hyperedge $\left\{u_1,v_1,v_n\right\}$ at time zero, as is the case in ${B'_0}^n$, then, for $0\leq t\leq n-2$, the copy of $K_4^3$ with this vertex set has at most two infected hyperedges at time $t$, and is therefore not susceptible at that time. Any other copy of $K_4^3$ containing $\left\{u_1,v_1,v_n\right\}$ has at most one infected hyperedge. Thus, we have that ${B'_t}^n=B_t^n\cup\left\{\left\{u_1,v_1,v_n\right\}\right\}$ for $0\leq t\leq n-1$. 

Given that the hyperedge $\left\{u_1,v_1,v_n\right\}$ is infected in ${B'_0}^n$, we see that the copy of $K_4^3$ with vertex set $\left\{u_1,u_2,v_1,v_n\right\}$ is susceptible to ${B'}_{n-1}^n$. Note that, by Lemma~\ref{lem:bb}, $\left\{u_1,v_1,v_n\right\}$ is the unique wide hyperedge (to borrow terminology from the proof of Lemma~\ref{lem:bb}) that is infected at time $n-1$. Thus, there is only one susceptible copy of $K_4^3$ at time $n-1$ and we have that ${B'_n}^n={B'}_{n-1}^n\cup \left\{\left\{u_2,v_1,v_n\right\}\right\}$. Let us show that ${B'_n}^n$ is $K_4^3$-stable. Note that the only wide or flat infected hyperedges at time $n$ are precisely $\left\{u_1,v_1,v_n\right\}$ and $\left\{u_2,v_1,v_n\right\}$, both of which are wide. Any $K_4^3$ containing both of $u_1$ and $u_2$ is either fully infected in ${B'_n}^n$ or contains two healthy wide hyperedges. Any $K_4^3$ containing exactly one of $u_1$ or $u_2$ contains at least one healthy wide hyperedge (here, we use that $n\geq4$) and exactly one healthy flat hyperedge. Finally, any $K_4^3$ consisting containing neither $u_1$ nor $u_2$ consists only of flat hyperedges and therefore none of its edges are in ${B'_n}^n$. Thus, ${B_0'}^n$ is $K_4^3$-tame. Note that the only $0\leq j\leq n-1$ for which $F_j$ contains the hyperedge $e_0'$ is $F_0$. Also, $\left\{u_1,u_2,v_1,v_n\right\}$ does not contain any of the hyperedges $e_0',e_1,\dots,e_{n-2}$. Combining this with Lemma~\ref{lem:bb}, we see that ${B'_0}^n$ satisfies condition \ref{civilized1} of Definition~\ref{defn:civilized} with respect to $e_0'$. Note that the only hyperedges that are present in ${B_0'}^n\setminus\left\{e_0'\right\}$ and absent from $B_0^n\setminus\left\{e_0\right\}$ are $e_0$ and $\left\{u_1,v_1,v_n\right\}$. Thus, since $B_0^n$ is $K_4^3$-civilized with respect to $e_0$, if there is a copy of $K_4^3$ that is susceptible to ${B_0'}^n\setminus\left\{e_0'\right\}$, then it must contain $e_0$ or $\left\{u_1,v_1,v_n\right\}$. However, every such copy has at least two hyperedges that are missing from ${B_0'}^n\setminus\left\{e_0'\right\}$, and so condition \ref{civilized2} of Definition~\ref{defn:civilized} is satisfied. 
\end{proof}

Next, we show that, by carefully chaining together augmented beachball hypergraphs, one obtains a $K_4^3$-civilized hypergraph with quadratic running time. This implies Theorem~\ref{th:n^2} via Corollary~\ref{cor:reduction}. 

\begin{theorem}
\label{th:K43}
For any $n\geq4$ and $m\geq 2$ there exists a $K_4^3$-civilized hypergraph $H_0$ with $n+m$ vertices such that
\[M_4^3\left(H_0\right)= \left(m-1\right)n.\]
\end{theorem}

\begin{proof}
Let $v_0,\dots,v_{n-1}$ be distinct vertices whose indices are viewed ``cyclically'' modulo $n$. In particular, for any $i$, we have $v_{-i}=v_{n-i}=v_{2n-i}$, and so on. Let $u_0,\dots,u_{m-1}$ be $m$ additional vertices which are distinct from one another, and from $v_0,\dots,v_{n-1}$. For $0\leq i\leq m-2$, define $S_i:=\left\{v_0,\dots,v_{n-1}\right\}\cup\left\{u_i,u_{i+1}\right\}$. We let $H_0$ be the hypergraph
\[B'\left(v_{0},v_{1},\dots,v_{n-1},u_0,u_{1}\right)\cup
\left(\bigcup_{i=1}^{m-2}B'\left(v_{-i},v_{1-i},\dots,v_{n-1-i},u_i,u_{i+1}\right)\setminus\left\{\left\{v_{-i},v_{1-i},u_i\right\}\right\}\right).\]
The construction can be thought of intuitively as follows. Start with the hypergraph $B'\left(v_0,\dots,v_{n-1},u_0,u_1\right)$. Then, we take the same beachball hypergraph again, except that we delete one specific hyperedge, rotate the middle vertices by one, move the bottom vertex to the top, and insert a new vertex on the bottom. Repeat until there are no remaining vertices to add. See Figure~\ref{fig:beachball2} for an illustration. 

\begin{figure}[htbp]
  \centering
    \includegraphics[width=1\textwidth]{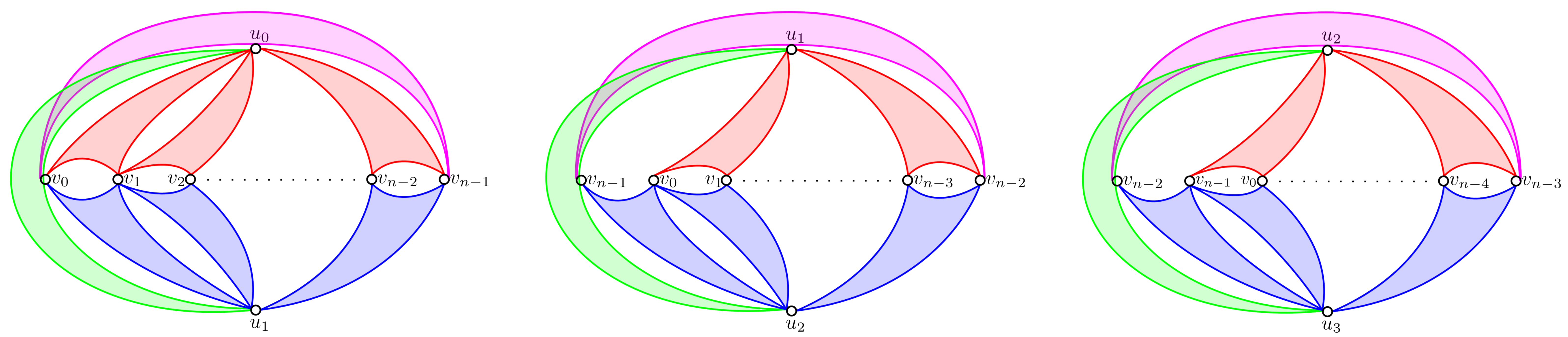}
    \caption{The hypergraph $H_{0}$ constructed in the proof of Theorem~\ref{th:K43}. We have drawn the only hyperedges of $H_0$ that are on $S_0$, $S_1$ and $S_2$ for illustration purposes.}
    \label{fig:beachball2}
\end{figure}

First, we claim that the hyperedges of $H_0$ within $S_0$ are precisely the hyperedges of the augmented beachball hypergraph on its vertices. Every hyperedge of the augmented beachball hypergraph on $S_0$ was added, by definition of $H_0$, and so it suffices to show that no additional hyperedges are present. For $i>1$, all of the hyperedges of  $B'\left(v_{-i},v_{1-i},\dots,v_{n-1-i},u_i,u_{i+1}\right)$ contain $u_i$ or $u_{i+1}$ and thus are not contained in $S_0$. Moreover, the only hyperedge of $B'\left(v_{n-1},v_{0},\dots,v_{n-2},u_1,u_{2}\right)$ contained in $S_0$ that is not a hyperedge of the augmented beachball on its vertices is $\left\{v_{n-1},v_0,u_1\right\}$ which was deleted during the construction of $H_0$. This proves the claim. The same argument also proves that, for all $i\geq1$, the hyperedges within the vertices $S_i$ are precisely those of the augmented beachball hypergraph on these vertices, except that the hyperedge $\left\{v_{-i},v_{1-i},u_i\right\}$ is missing. 

Now, we also observe that, for any $t\geq0$, the hypergraph $H_t$ does not contain any hyperedge which includes $u_i$ and $u_{i'}$ with $|i-i'|>1$. To see this, consider the first time that it fails and consider the susceptible $K_4^3$ which caused it; as $H_0$ contains no such hyperedges containing $u_i$ and $u_{i'}$, we get a contradiction. 

By Lemma~\ref{lem:augbb}, every augmented beachball in the chain is $K_4^3$-civilised, and the corresponding trajectory implies that the hyperedge $\left\{v_{-i},v_{1-i},u_i\right\}$ is contained in the unique $K_4^3$ that is susceptible to $B'\left(v_{-i},v_{1-i},\dots,v_{n-1-i},u_i,u_{i+1}\right)$. 

So we see that, for any $i\geq1$, the first hyperedge in $S_i$ that becomes infected does so due to a susceptible copy of $K_4^3$ with precisely one vertex outside of $S_i$ (since it cannot be contained entirely inside $S_i$); until such a copy appears, no additional hyperedges in $S_i$ become infected. Note that this statement is not true for $i=0$, as the first hyperedge in $S_0$ becomes infected due to a susceptible copy of $K_4^3$ contained within $S_0$.

Putting this together, we see that the only hyperedges that become infected in the first $n$ steps are exactly those which become infected in the $K_4^3$-bootstrap process in $B'\left(v_{0},v_{1},\dots,v_{n-1},u_0,u_{1}\right)$. By Lemma~\ref{lem:augbb}, all of these hyperedges are disjoint from $S_i$ for all $i\geq1$, except for the last one, namely $\left\{u_1,v_0,v_{n-1}\right\}$, which happens to be the unique hyperedge of $B'\left(v_{n-1},v_{0},\dots,v_{n-2},u_1,u_{2}\right)$ that is missing from $H_0$. From this point, the infection follows the $K_4^3$-bootstrap process in $B'\left(v_{n-1},v_{0},\dots,v_{n-2},u_1,u_{2}\right)$ for $n$ steps. This pattern repeats itself $m-1$ times. Thus, $H_0$ is $K_4^3$-tame. The proof that the additional conditions of Definition~\ref{defn:civilized} hold with respect to $e_0'=\left\{u_0,v_0,v_1\right\}$ is analogous to the arguments given in the proofs of the lemmas in this section, and so we omit it. This completes the proof.
\end{proof}

\begin{proof}[Proof of Theorem~\ref{th:n^2}]
For each $n\geq8$, Theorem~\ref{th:K43} gives us a $K_4^3$-civilized hypergraph $H_0$ on $n$ vertices such that
\[M_4^3\left(H_0\right)\geq \left(\left\lceil\frac{n}{2}\right\rceil-1\right)\left\lfloor\frac{n}{2}\right\rfloor.\]
The result now follows by Corollary~\ref{cor:reduction}.
\end{proof}

\section{Cubic Time Construction for \texorpdfstring{$\boldsymbol{K_5^3}$}{K 5 3}}
\label{sec:cubic}

Our final task is to use the beachball hypergraph again to obtain a construction for the $K_5^3$-bootstrap process with cubic running time. 

\begin{theorem}
\label{th:K53}
For any $n\geq2$ and $m\geq2$, there exists a $K_5^3$-civilized hypergraph $H_0$ with $2n+4m+4$ vertices such that
\[M_5^3\left(H_0\right)=\left(2n+4\right)m^2 - 1.\]
\end{theorem}

\begin{proof}
The vertices of the construction are naturally divided into two ``halves'' which share exactly two special vertices, $w_1$ and $w_2$. For $i\in \left\{1,2\right\}$, the vertices that are ``exclusive'' to the $i$th half are $v_1^i,\dots, v_n^i,u_1^i,\dots,u_m^i,z_1^i,\dots,z_m^i$ and $x^i$. Each half of the construction has $n+2m+1$ exclusive vertices. Together with $w_1$ and $w_2$, this makes $2n+4m+4$ vertices in total. 

For $1\leq j\leq m^2$, let $\left(a_j,b_j\right)$ be the $j$th element of $\left\{1,2,\dots,m\right\}\times \left\{1,2,\dots,m\right\}$ in lexicographic order. For $1\leq j\leq m^2$, let $C_j^1$ be the copy of $K_5^3$ with vertex set $\left\{u_{a_j}^1,z_{b_j}^1,u_{a_j}^2,z_{b_j}^2,w_2\right\}$ and, if $j\leq m^2-1$, then we additionally define $C_j^2$ to be the copy of $K_5^3$ with vertex set $\left\{u_{a_j}^2,z_{b_j}^2,u_{a_{j+1}}^1,z_{b_{j+1}}^1,w_1\right\}$. Let 
\[B^{1}_j := B\left(w_1,v_1^1,\dots,v_n^1,w_2,u_{a_j}^1,z_{b_j}^1\right)\vee x^1\]
and 
\[B^{2}_j:=B\left(w_2,v_1^2,\dots,v_n^2,w_1,u_{a_j}^2,z_{b_j}^2\right)\vee x^2\]
for all $1 \leq j \leq m^2$. Define $H_0$ to be the union of the following five hypergraphs
\[B^{1}_1,\]
\[\bigcup_{j=2}^{m^2}B^{1}_j\setminus\left\{\left\{u_{a_j}^1,z_{b_j}^1,w_1\right\}\right\},\]
\[\bigcup_{j=1}^{m^2}B^{2}_j\setminus\left\{\left\{u_{a_j}^2,z_{b_j}^2,w_2\right\}\right\},\]
\[\bigcup_{j=1}^{m^2}C_j^1\setminus\left\{\left\{u_{a_j}^1,z_{b_j}^1,w_2\right\},\left\{u_{a_j}^2,z_{b_j}^2,w_2\right\}\right\}, \text{ and}\]
\[\bigcup_{j=1}^{m^2-1}C_j^{2}\setminus\left\{\left\{u_{a_j}^2,z_{b_j}^2,w_1\right\},\left\{u_{a_{j+1}}^1,z_{b_{j+1}}^1,w_1\right\}\right\}.\]

Intuitively, the way that the infection evolves can be described as follows. It starts by propagating through $B_1^1$, the first beachball hypergraph on the first half, in the way described by Lemma~\ref{lem:bb} and Proposition~\ref{prop:easyStepUp}. This takes $n+1$ steps. The hyperedge that gets infected in the last step is $\left\{u_{a_1}^1,z_{b_1}^1,w_2\right\}$. Consequently, now all of the hyperedges within $V\left(C_1^1\right)$ except for $\left\{u_{a_1}^2,z_{b_1}^2,w_2\right\}$ are infected, and so this hyperedge gets infected in the next step; at this point, we have done $n+2$ steps in total. This was the only hyperedge of $B_1^2$ that is not in $H_0$, and so its infection triggers the beachball hypergraph $B_1^2$ on the right half of the construction. The process inside $B_1^2$ ends with the infection of the hyperedge $\left\{u_{a_1}^2,z_{b_1}^2,w_1\right\}$, which is in $V\left(C_1^2\right)$. So now all the hyperedges of $V\left(C_1^2\right)$ except $\left\{u_{a_2}^1,z_{b_2}^1,w_1\right\}$ are infected, and hence this hyperedge gets infected in the next step. This transfers the infection back over to the left half of the construction by triggering $B_2^1$, and so on. If the process does indeed progress in this manner, then, in total, it will take $n+2$ steps for every $1\leq j\leq m^2$ and every $i\in\{1,2\}$, with the exception of $j=m^2$ and $i=2$, which only contributes $n+1$ steps (since $C_{m^2}^2$ has not been defined). Thus, there are $\left(2n+4\right)m^2-1$ steps. 

Let us now make this rigorous. Our goal is to show that $H_0$ is $K_5^3$-civilized with respect to $e_0=\left\{u_{a_1}^1,z_{b_1}^1,w_1\right\}$; on the way to that goal, we will need to establish several claims. First, we observe that any pair of vertices which are not contained together in a hyperedge of $H_0$ are not contained together in any hyperedge of $H_t$ for any $t\geq0$; to see this, consider the first time that such a hyperedge becomes infected, look at the susceptible copy of $K_5^3$ which caused it and get a contradiction. This argument immediately yields the following claim. Throughout the statement of the claim, keep in mind that $w_1$ and $w_2$ are regarded as being on both halves of the construction, but all other vertices are exclusive to one half or the other.

\begin{claim}
\label{claim:badEdges1}
For any $t\geq0$, $H_t$ does not contain any hyperedge $e$ which satisfies any of the following conditions:
\begin{enumerate}
    \item[(i)]\label{xGood} $e$ contains $x^i$ but is not contained in the $i$th half of the construction for some $i\in \left\{1,2\right\}$,
    \item[(ii)]\label{vGood} $e$ contains $v_j^i$ but is not contained in the $i$th half of the construction for some $1\leq j\leq n$ and  $i\in \left\{1,2\right\}$,
    \item[(iii)]\label{uzGood} $e$ contains two of $u_1^i,\dots,u_m^i$ or two of $z_1^i,\dots,z_m^i$ for some $i\in\left\{1,2\right\}$.
\end{enumerate}
\end{claim}

The purpose of the next claim is similar to the previous one; i.e. to rule out certain types of hyperedges from becoming infected.

\begin{claim}
\label{claim:badEdges2}
For any $t\geq0$, every hyperedge of $H_t$ containing both $w_1$ and $w_2$ must also contain one of $x^1$ or $x^2$.
\end{claim}

\begin{proof}[Proof of Claim~\ref{claim:badEdges2}]
Consider the first time $t$ such that $H_t$ contains a hyperedge $e$ that contains both of $w_1$ and $w_2$ but neither of $x^1$ nor $x^2$. Note that $t\geq1$ by construction of $H_0$. Let $F$ be the copy of $K_5^3$ containing $e$ which is susceptible to $H_{t-1}$. By minimality of $t$, the only two hyperedges of $H_{t-1}$ containing $w_1$ and $w_2$ are $\left\{w_1,w_2,x^1\right\}$ and $\left\{w_1,w_2,x^2\right\}$. So, the vertices of $F$ must consist of $w_1,w_2,x^1,x^2$ and a fifth vertex, say $y$. However, $y$ is exclusive to one of the two halves of the construction, and so, regardless of which half it is, $H_{t-1}$ contains several hyperedges which satisfy condition \ref{xGood} of Claim~\ref{claim:badEdges1}, which is a contradiction.
\end{proof}

Let us now show that $H_0$ is $K_5^3$-tame with the trajectory that was described earlier in the proof. Suppose that this is not the case, let $t\geq1$ be the minimum time in which there is an unexpected infected hyperedge, say $e$. Let $F$ be a copy of $K_5^3$ containing $e$ which is susceptible to $H_{t-1}$. By Lemma~\ref{lem:bb} and minimality of $t$, we cannot have $V\left(F\right)\subseteq B_j^i$ for any $1\leq j\leq m^2$ and $i\in\{1,2\}$. Now, any five vertices on the $i$th half of the construction are either contained in $B_j^i$ for some $j$ or contain two of $u_1^i,\dots,u_m^i$ or two of $z_1^i,\dots,z_m^i$. So, as there are no hyperedges satisfying condition \ref{uzGood} of Claim~\ref{claim:badEdges1}, $V\left(F\right)$ must contain at least one vertex that is exclusive to each half of the construction. But now, since no hyperedges satisfy Claim~\ref{claim:badEdges1} \ref{xGood} or \ref{vGood}, we get that $V\left(F\right)$ cannot contain any of the vertices $x^i$ or $v_j^i$ for $i\in\left\{1,2\right\}$ and $1\leq j\leq n$. This also means that it cannot contain both of $w_1$ and $w_2$ by Claim~\ref{claim:badEdges2}. On the other hand, it must contain at least one of $w_1$ or $w_2$; if not, then three vertices of $F$ are exclusive to one of the sides of the partition, and we get a contradiction from Claim~\ref{claim:badEdges1} \ref{uzGood}. Using Claim~\ref{claim:badEdges1} \ref{uzGood} one more time, we see that $V\left(F\right)$ has the form $\left\{u_a^1,z_b^1,w_\ell,u_c^2,z_d^2\right\}$ for some $1\leq a,b,c,d\leq m$ and $\ell\in\left\{1,2\right\}$. We now divide the proof into cases. 

\begin{case}
$e$ does not contain $w_\ell$. 
\end{case}

We assume that $e=\left\{u_a^1,z_b^1,u_c^2\right\}$ and note that the other three cases follow from similar arguments. Let $1\leq j\leq m^2$ be chosen so that $a=a_j$ and $b=b_j$; such a $j$ exists by construction of $H_0$. The hypergraph $H_{t-1}$ contains the hyperedge $e'=e\setminus\left\{u_c^2\right\}\cup \left\{z_d^2\right\}$ and so, by minimality of $t$, we get that $e'$ is contained in either $V\left(C_j^1\right)$ or $V\left(C_{j-1}^2\right)$. Indeed, under the trajectory described at the beginning of the proof, every hyperedge containing two vertices exclusive to one side and one exclusive to the other had this property. The fact that $e'$ is contained in $V\left(C_j^1\right)$ or $V\left(C_{j-1}^2\right)$ implies that $e$ is  as well. However, every hyperedge of $C_j^1$ or $C_{j-1}^2$ not containing $w_1$ or $w_2$ was added to $H_0$ originally, which contradicts our choice of $e$.

\begin{case}
$e$ contains $w_\ell$.
\end{case}

Again, let $1\leq j\leq m^2$ be chosen so that $a=a_j$ and $b=b_j$. The hyperedge $\left\{u_{a_j}^1,z_{b_j}^1,u_c^2\right\}$ is contained in $H_{t-1}$ and so, by minimality of $t$, it must be contained in either $V\left(C_j^1\right)$ or $V\left(C_{j-1}^2\right)$. By definition of $V\left(C_j^1\right)$ and $V\left(C_{j-1}^2\right)$, and our specific choice of lexicographic order in the construction of $H_0$, this implies that either $c=a_j$ or $c=a_j-1$ and $b_j=1$. 

Suppose first that $c=a_j$. Then we get that $H_{t-1}$ contains $\left\{z_{b_j}^1,u_{a_j}^2,z_d^2\right\}$ which, by minimality of $t$ and construction of $H_0$, implies that either $d=b_j$ or $d=b_j-1$. If $d=b_j$, then we must have that $\ell=2$ by minimality of $t$ since $\left\{z_{b_j}^1,w_1,z_{b_j}^2\right\}$ is not contained in any of the hypergraphs $C_1^1,\dots,C_{m^2}^1$ or $C_1^2,\dots,C_{m^2-1}^2$. So, in the case that $c=a_j$ and $d=b_j$, we get that $V\left(F\right)=V\left(C_j^1\right)$. This implies that $e$ is one of the only two hyperedges of $C_j$ that are missing from $H_0$. If $e=\left\{u_{a_j}^1,z_{b_j}^1,w_2\right\}$, then the hyperedge $\left\{w_2,u_{a_j}^2,z_{b_j}^2\right\}$ is present before $e$, which contradicts minimality of $t$. If $e=\left\{w_2,u_{a_j}^2,z_{b_j}^2\right\}$, then $e$ is being infected due to $C_j^1$ being susceptible, which fits the description of the trajectory from earlier in the proof, and so this contradicts the definition of $e$. Now, suppose that $c=a_j$ and $d=b_j-1$. In this case, we must have $\ell=1$ by minimality of $t$ since $\left\{z_{b_j}^1,w_2,z_{b_j-1}^2\right\}$ is not contained in any of the sets $C_1^1,\dots,C_{m^2}^1$ or $C_1^2,\dots,C_{m^2-1}^2$. So, what we end up with is $V\left(F\right)=V\left(C_{j-1}^2\right)$ and we get a contradiction similar to the case that we just analyzed.

Now, suppose that $c=a_{j}-1$ and $b_j=1$. Since $H_{t-1}$ contains $\left\{z_{1}^1,u_{a_{j}-1}^2,z_d^2\right\}$, we must have either $d=1$ or $d=m$. If $d=1$, we get a contradiction, since $\left\{u_{a_j}^1,u_{a_j-1}^2,z_1^2\right\}$ is not contained in $H_{t-1}$ by minimality of $t$ and the fact that $m\geq2$ and so, under lexicographic order, $\left(a_j-1,1\right)$ does not immediately precede or follow any pair involving $a_j$ in the first coordinate. So, $d=m$. Now, we get that $\ell=1$ by minimality of $t$ since $\left\{z_1^1,w_2,z_m^2\right\}$ is not contained in any of the sets $C_1^1,\dots,C_{m^2}^1$ or $C^2_1,\dots,C^2_{m^2-1}$. So, we get that $V\left(F\right)= C^2_{j-1}$ which leads us to a contradiction, as in the previous paragraph. 

Similar arguments also show that, for $e_0=\left\{u_{a_1}^1,z_{a_1}^1,w_1\right\}$, the hypergraph $H_0\setminus\left\{e_0\right\}$ is $K_5^3$-stable; thus, condition \ref{civilized2} of Defintion~\ref{defn:civilized} holds. Clearly each of the copies of $K_5^3$ in the trajectory of $H_0$, as described above, has precisely two hyperedges missing from $H_0\setminus\left\{e_0\right\}$ and so \ref{civilized1} holds, too. This completes the proof.
\end{proof}

\begin{proof}[Proof of Theorem~\ref{th:n^3}]
For $n$ sufficiently large and $n\equiv 4\bmod 6$, Theorem~\ref{th:K53} provides a $K_5^3$-civilized hypergraph $H_0'$ with $n$ vertices such that
\[M_5^3\left(H_0'\right)\geq \left(\frac{n+8}{3}\right)\left(\frac{n-4}{6}\right)^2-1.\]
The result follows by Corollary~\ref{cor:reduction}.
\end{proof}

\section{Open Problems}
\label{sec:concl}

While Theorem~\ref{th:n^3} determines $M_k^r\left(n\right)$ up to a constant factor for all $r\geq 3$ and $k\geq r+2$, it would be interesting to have a better understanding of the growth rate in the case $k=r+1$. The most interesting question here seems to be whether there is a non-trivial upper bound on the running time of the $K_4^3$-bootstrap process; we conjecture that the quadratic lower bound from Theorem~\ref{th:n^2} is tight up to a constant factor.

\begin{conjecture}
\label{conj:K43}
$M_4^3\left(n\right)=O\left(n^2\right)$. 
\end{conjecture}

In contrast, we do not believe that Theorem~\ref{th:n^2} is tight for all $r\geq3$. It is even conceivable that, for large enough $r$, the maximum running time of the $K_{r+1}^r$-process on $n$ vertices eventually becomes $\Theta\left(n^r\right)$. We ask whether this is, indeed, the case.

\begin{question}
\label{ques:gap1}
Does there exist an integer $r_0$ such that, if $r\geq r_0$, then $M_{r+1}^r\left(n\right)=\Theta\left(n^r\right)$? 
\end{question}

\section{Acknowledgements}
We would like to thank the anonymous referee for their helpful comments that improved particular details as well as the overall coherence of the paper.

\section*{Remarks}
After completing this work, Hartarsky and Lichev~\cite{HartarskyLichev22+} and Espuny D\'{i}az, Janzer, Kronenberg and Lada~\cite{Espuny+22+} independently disproved~\ref{conj:K43} by showing that $M_4^3(n) = \Theta(n^3)$, and consequently answered Question~\ref{ques:gap1} in the affirmative in a rather strong sense with $r_0=3$. Hartarsky and Lichev~\cite{HartarskyLichev22+} also determine the leading asymptotics of the prefactor when $r \to \infty$. Additionally, Espuny D\'{i}az et al.~\cite{Espuny+22+} provide the first nontrivial exact result about the maximum running times of hypergraph bootstrap percolation by showing that the maximum running time for the $H$-bootstrap process on $n$ vertices when $H$ is $K_4^3$ minus an edge is exactly $2n- \lfloor \log_2(n-2) \rfloor -6$.

\newblock 

\bibliographystyle{plain}

\end{document}